\documentclass[a4paper,11pt,twoside,reqno]{amsart}

\usepackage[utf8]{inputenc}
\usepackage[plainpages=false,pdfpagelabels=true]{hyperref}
\usepackage{amssymb,amsthm}
\usepackage[margin=1in]{geometry}

\newtheorem{Theorem}{Theorem}[section]

\newtheorem{Lem}[Theorem]{Lemma}

\theoremstyle{definition}

\newtheorem{Bem}[Theorem]{Remark}

\newcommand{\tr}{\operatorname{Tr}}
\newcommand{\Int}{\operatorname{Int}}

\newcommand{\dv}{\text{ }dV}
\allowdisplaybreaks[1]

\renewcommand{\epsilon}{\varepsilon}

\newcommand{\R}{\ensuremath{\mathbb{R}}}

\numberwithin{equation}{section}

\title{Unique continuation theorems for biharmonic maps}

\author{Volker Branding}

\address{University of Vienna, Faculty of Mathematics\\
Oskar-Morgenstern-Platz 1, 1090 Vienna, Austria}
\email{volker.branding@univie.ac.at}

\author{Cezar Oniciuc}

\address{Faculty of Mathematics\\
Al.I. Cuza University of Iasi, Bd. Carol I, 11, 700506 Iasi, Romania}
\email{oniciucc@uaic.ro}

\date{\today}

\subjclass[2010]{58E20; 31B30}

\keywords{biharmonic maps; unique continuation}

\begin{document}

\begin{abstract}
We prove several unique continuation results for biharmonic
maps between Riemannian manifolds.
\end{abstract}

\maketitle


\section{Introduction and results}
Finding interesting maps between Riemannian manifolds is one of the most
challenging problems in modern Riemannian geometry.
Suppose that \((M,g)\) and \((N,h)\) are two Riemannian manifolds.
Moreover, let \(\phi\colon M\to N\) be a smooth map.
One option of finding such maps is to find extrema of their energy
\begin{align}
\label{energy-harmonic}
E(\phi)=\int_M\left|d\phi\right|^2\dv_g,
\end{align}
which are called \emph{harmonic maps}
and are characterized by the vanishing of the so-called \emph{tension field}, that is
\begin{align}
\label{harmonic}
0=\tau(\phi):=\tr_g\nabla d\phi,
\end{align}
where \(\tau(\phi)\in\Gamma(\phi^\ast TN)\).
Many results on harmonic maps have been obtained in the past decades,
we refer to \cite{MR2389639} for an overview.

Currently, there is a growing interest in a geometric variational problem
that generalizes harmonic maps, the so-called \emph{biharmonic maps},
which were first studied in \cite{MR886529}.
In this case one looks for critical points of the bienergy of a map \(\phi\colon M\to N\),
which is defined as
\begin{align}
\label{energy-biharmonic}
E_2(\phi)=\int_M\left|\tau(\phi)\right|^2\dv_g.
\end{align}
The critical points of \eqref{energy-biharmonic} are characterized
by the vanishing of the so-called \emph{bitension field}, that is
\begin{align}
\label{biharmonic}
0=\tau_2(\phi):=\Delta\tau(\phi)-\sum_{i=1}^mR^N\left(d\phi\left(e_i\right),\tau(\phi)\right)d\phi\left(e_i\right),
\end{align}
where \(\{e_i\},i=1,\ldots,m=\dim M\) is a local orthonormal frame field tangent to \(M\) and \(\Delta\) represents
the Laplacian on \(\phi^\ast TN\).
For recent results on biharmonic submanifolds we refer to \cite{cezar-habil}, see also the older survey \cite{MR2301373}.
The harmonic map equation \eqref{harmonic} is a second order elliptic
partial differential equation, whereas the
biharmonic map equation \eqref{biharmonic} is of fourth order,
which makes it substantially harder to characterize the
qualitative behavior of its solutions.

In this article we want to focus on one particular aspect regarding
the qualitative behavior of harmonic and biharmonic maps, namely the unique continuation property.
For harmonic maps the question of unique continuation was
settled by Sampson \cite[Theorem 1]{MR510549}, see also \cite[Section 1.4.2]{MR1391729}.
More precisely, the following result was proved:

\begin{Theorem}[Sampson]
\label{sampson}
Let \(\phi_1,\phi_2\colon M\to N\) be two harmonic maps.
If they agree on an open subset
then they are identical; and indeed the conclusion
holds if \(\phi_1\) and \(\phi_2\) agree to infinitely high order at some point.
In particular, a harmonic map which is constant on an open subset is a constant map.
\end{Theorem}

In addition, Sampson established the following geometric unique continuation
property of harmonic maps (see \cite[Theorem 6]{MR510549}):
\begin{Theorem}[Sampson]
\label{sampson-submanifold}
Let \(\phi\colon M\to N\) be a harmonic map and let \(P\) be a regular, closed,
totally geodesic submanifold of \(N\). If an open set of \(M\) is mapped into \(P\),
then all of \(M\) is mapped into \(P\).
\end{Theorem}

There is a strong belief that the unique continuation property holds
for all solutions of second order elliptic partial differential equations
arising in geometry \cite{MR948075}. On the other hand
one can construct explicit counterexamples of solutions of fourth order elliptic
partial differential equations where the unique continuation
property does not hold \cite[Example 1.11]{MR948075}.

In this article we will prove the following results for biharmonic maps:

\begin{Theorem}
\label{harmonic-everywhere}
Let $\phi\colon M\to N$ be a biharmonic map. If $\phi$ is harmonic on an open subset, then it is harmonic everywhere.
\end{Theorem}

Since there is a big interest in biharmonic maps to spheres we will first establish a unique continuation
result for spherical targets as also the proof is considerably simpler here.

\begin{Theorem}
\label{sphere}
Let \(\phi_1,\phi_2\colon M\to \mathbb{S}^n\) be two biharmonic maps. If they agree on an open subset, then they are identical.
\end{Theorem}

Afterwards, we will prove a unique continuation theorem for biharmonic maps to an arbitrary target manifold.

\begin{Theorem}
\label{general-target}
Let \(\phi_1,\phi_2\colon M\to N\) be two biharmonic maps. If they agree on an open subset, then they are identical.
\end{Theorem}

Finally, we will also provide a geometric unique continuation property for biharmonic maps
generalizing Sampson's result for harmonic maps \cite[Theorem 6]{MR510549}.
First, we will give the following version for a spherical target:

\begin{Theorem}
\label{theorem-totally-geodesic-sphere}
Let \(\phi\colon M\to \mathbb{S}^n\) be a biharmonic map.
If an open subset of \(M\) is mapped into the equator \(\mathbb{S}^{n-1}\),
then all of \(M\) is mapped into \(\mathbb{S}^{n-1}\).
\end{Theorem}

In addition, we are also able to prove a corresponding version of the above theorem
for an arbitrary target.
\begin{Theorem}
\label{theorem-totally-geodesic}
Let \(\phi\colon M\to N\) be a biharmonic map and let \(P\) be a regular, closed, totally geodesic
submanifold of \(N\). If an open subset of \(M\) is mapped into \(P\) then all of \(M\)
is mapped into \(P\).
\end{Theorem}

Finally, we show by an explicit counterexample that Sampson's maximum principle for harmonic maps
\cite[Theorem 2]{MR510549} does not extend to biharmonic maps.

\begin{Bem}
\begin{enumerate}
\item Theorems \ref{sphere} and \ref{theorem-totally-geodesic} were proved in \cite[Theorem 4]{MR3523535} and \cite[Proposition 3]{MR3523535}
for the particular case of CMC biharmonic immersions into spheres. In the given reference the following idea was employed:
A CMC biharmonic immersion into \(\mathbb{S}^n\) composed with the canonical inclusion of the sphere in the ambient Euclidean space \(\R^{n+1}\)
gives an immersion that can be written as a sum of two \(\R^{n+1}\)-valued eigenmaps of the Laplace operator on \((M^m,g)\).
These maps induce harmonic maps into \(\mathbb{S}^n\) of appropriate radius. Then, in contrast to the article at hand, the results
follow directly from the classical results of Sampson \cite[Theorems 1,6]{MR510549} for harmonic maps.
Here, our results are more general and the technique we are using is completely different.
\item Theorem \ref{harmonic-everywhere} was first proved in \cite{MR1919374} by simply applying \cite[Proposition 1.2.3]{MR0233295},
but here, the proof is more clear and based on the classical result from Aronszajn \cite{MR0092067}.
\item In \cite[Theorem 5.3]{MR2332421} a unique continuation result for extrinsic biharmonic
maps from \(\Omega\subset\R^4\) to \(S^4\) was proved.
\item Theorems \ref{sphere} - \ref{theorem-totally-geodesic} also hold for extrinsic biharmonic maps and semi-biharmonic maps, which were introduced in \cite{branding2018semi}.
\end{enumerate}
\end{Bem}

Throughout this article we will use the following sign conventions.
For the Riemannian curvature
tensor field we use \(R(X,Y)Z=[\nabla_X,\nabla_Y]Z-\nabla_{[X,Y]}Z\) and
for the (rough) Laplacian on \(\phi^\ast TN\) we use
\(\Delta:=\tr_g(\nabla\nabla-\nabla_\nabla)\).

In general, we will use the same symbol $\langle\cdot,\cdot\rangle$ to indicate the Riemannian metrics on various vector bundles, and the same symbol $\nabla$ for the corresponding Riemannian connections.

All manifolds are assumed to be connected and we will work only with smooth objects.

Whenever we will make use of indices, we will use
Latin indices \(i,j,k\) for indices on the domain ranging from \(1\) to \(m\)
and Greek indices \(\alpha,\beta,\gamma\) for indices on the target
which take values between \(1\) and \(n\). When the range of the indices is from $1$ to $q$, for some positive integer $q$, we will often denote them
by \(a,b,c\).

We will use the Einstein summation convention, i.e. repeated indices that are in the diagonal position indicate the sum.

\section{Proofs of the theorems}

In this section we will prove the results obtained in this article.

Our strategy of proof is to cleverly rewrite the biharmonic map
equation such that we are effectively dealing with a second order problem
to which we can apply the classical result from Aronszajn \cite{MR0092067}
extending the ideas from \cite{MR2332421}.

\subsection{Proof of Theorem \ref{harmonic-everywhere}}

We recall the following (\cite[p.248]{MR0092067})

\begin{Theorem}
\label{aro-theorem}
Let \(A\) be a linear elliptic second-order differential operator defined on an open subset \(D\) of \(\R^m\).
Let \(u=(u^1,\ldots,u^q)\) be functions in \(D\) satisfying the inequality
\begin{equation}
\label{aro-voraus}
\left|Au^a\right|\leq C\left(\sum_{b,i}\left|\frac{\partial u^{b}}{\partial x^i}\right|+\sum_{b}\left|u^{b}\right|\right).
\end{equation}
If \(u=0\) in an open subset of $D$, then \(u=0\) throughout \(D\).
\end{Theorem}

Let $\phi:M\to N$ be a smooth map and let $\sigma$ be a section in the pull-back bundle $\phi^{\ast}TN$.
We consider a local chart $(U,x^i)$ on $M$ and a local chart $(V,y^{\alpha})$ on $N$ such that $\phi(U)\subset V$. The section $\sigma$ can be written as
$$
\sigma = u^{\alpha}(p)\frac{\partial}{\partial y^{\alpha}}(\phi(p)), \quad p\in U.
$$
The Laplacian $\Delta\sigma$ is given by
\begin{align*}
\Delta\sigma = & \tr_g\nabla d\sigma = g^{ij}\left( \nabla d\sigma \right) \left( \frac{\partial}{\partial x^i},\frac{\partial}{\partial x^j}\right) \\
= & g^{ij}\left( \nabla_{\frac{\partial}{\partial x^i}}d\sigma\left(\frac{\partial}{\partial x^j}\right) - d\sigma\left( \nabla_{\frac{\partial}{\partial x^i}}\frac{\partial}{\partial x^j} \right) \right).
\end{align*}
We find
$$
d\sigma\left( \frac{\partial}{\partial x^i} \right) = \nabla_{\frac{\partial}{\partial x^i}}\sigma =  u^{\alpha}_i\frac{\partial}{\partial y^{\alpha}} + u^{\alpha} \phi^{\beta}_i \Gamma^{\theta}_{\beta \alpha}\frac{\partial}{\partial y^{\theta}},
$$
where $u^{\alpha}_i=\frac{\partial u^{\alpha}}{\partial x^i}$ and $\left(\phi^{\beta}\right)_\beta$ is the corresponding expression for $\phi$ in local coordinates.
Then by a straightforward computation we get (see also \cite[Lemma 1.1]{MR3011327})
\begin{align*}
\Delta\sigma - \tr_g R^N(d\phi,\sigma)d\phi = & \left\{ \Delta u^{\theta} + 2g^{ij} u^{\alpha}_j \phi^{\beta}_i \Gamma^{\theta}_{\beta \alpha} \right. \\
& \ + u^{\alpha}\left( \left( \Delta \phi^{\beta} \right) \Gamma^{\theta}_{\beta \alpha} \right.\\
& \ \ \ \ \ \ \left. \left. + g^{ij} \phi^{\beta}_j \phi^{\omega}_i \left( \frac{\partial \Gamma^{\theta}_{\beta \alpha}}{\partial y^{\omega}} + \Gamma^{\gamma}_{\beta \alpha}\Gamma^{\theta}_{\omega \gamma} - R^{\theta}_{\omega \beta \alpha} \right) \right) \right\}\frac{\partial}{\partial y^{\theta}}.
\end{align*}
We can go a bit further and notice that
\begin{align}
\label{local-coordinates}
g^{ij} \phi^{\beta}_j \phi^{\omega}_i &\left( \frac{\partial \Gamma^{\theta}_{\beta \alpha}}{\partial y^{\omega}} + \Gamma^{\gamma}_{\beta \alpha}\Gamma^{\theta}_{\omega \gamma}
- R^{\theta}_{\omega \beta \alpha} \right) \frac{\partial}{\partial y^{\theta}} \\
\nonumber= & \frac{1}{2} \left\langle d \phi^{\beta},d \phi^{\omega}\right\rangle \left( \frac{\partial \Gamma^{\theta}_{\beta \alpha}}{\partial y^{\omega}}
+ \Gamma^{\gamma}_{\beta \alpha}\Gamma^{\theta}_{\omega \gamma} - R^{\theta}_{\omega \beta \alpha} \right.
 \left. + \frac{\partial \Gamma^{\theta}_{\omega \alpha}}{\partial y^{\beta}} + \Gamma^{\gamma}_{\omega \alpha}\Gamma^{\theta}_{\beta \gamma} - R^{\theta}_{\beta \omega \alpha} \right) \frac{\partial}{\partial y^{\theta}} \\
\nonumber= & \frac{1}{2} \left\langle d \phi^{\beta},d \phi^{\omega}\right\rangle \left( \nabla_{\frac{\partial}{\partial y^{\alpha}}}\nabla_{\frac{\partial}{\partial y^{\beta}}}\frac{\partial}{\partial y^{\omega}} \right.
\left. + \nabla_{\frac{\partial}{\partial y^{\alpha}}}\nabla_{\frac{\partial}{\partial y^{\omega}}}\frac{\partial}{\partial y^{\beta}} \right) \\
\nonumber= & \left\langle d \phi^{\beta},d \phi^{\omega}\right\rangle \nabla_{\frac{\partial}{\partial y^{\alpha}}}\nabla_{\frac{\partial}{\partial y^{\beta}}}\frac{\partial}{\partial y^{\omega}} \\
\nonumber= & \left\langle d \phi^{\beta},d \phi^{\omega}\right\rangle \left( \frac{\partial \Gamma^{\theta}_{\beta \omega}}{\partial y^{\alpha}} + \Gamma^{\gamma}_{\beta \omega}\Gamma^{\theta}_{\alpha \gamma} \right) \frac{\partial}{\partial y^{\theta}}.
\end{align}
Thus, we obtain the following

\begin{Lem}
For $\phi\colon M\to N$ and $\sigma = u^{\alpha}\frac{\partial}{\partial y^{\alpha}}$ a section in $\phi^{\ast}TN$, we have
\begin{equation}
\label{local-coordinates-lemma}
\begin{aligned}
\Delta\sigma - \tr_g R^N(d\phi,\sigma)d\phi = & \left\{ \Delta u^{\theta} + 2\left\langle d  u^{\alpha},d  \phi^{\beta}\right\rangle \Gamma^{\theta}_{\alpha \beta} \right. \\
& \ + u^{\alpha}\left( \left( \Delta \phi^{\beta} \right) \Gamma^{\theta}_{\alpha \beta} \right.\\
& \ \ \ \ \ \ \ \ \ \left. \left. + \left\langle d \phi^{\beta},d \phi^{\omega}\right\rangle \left( \frac{\partial \Gamma^{\theta}_{\beta \omega}}{\partial y^{\alpha}} + \Gamma^{\gamma}_{\beta \omega}\Gamma^{\theta}_{\alpha \gamma} \right) \right) \right\}\frac{\partial}{\partial y^{\theta}}.
\end{aligned}
\end{equation}
\end{Lem}

\begin{proof}[Proof of Theorem \ref{harmonic-everywhere}]
Let us denote
$$
A:=\{p\in M : \tau(\phi)(p)=0\}.
$$
Clearly, $A$ is closed in $M$ and its topological interior, $\Int \ A$, is non-empty.

If the boundary of $\Int \ A$, $\partial (\Int \ A)$, is empty, then, as $M$ is connected, $\Int \ A = A = M$, i.e. $\phi$ is harmonic everywhere.

Assume that there exists $p_0\in \partial(\Int \ A)$. Furthermore, let $U$ be an arbitrary open subset containing $p_0$. Clearly $p_0$ does not belong to $\Int A$ and $U\cap \Int A \neq \emptyset$.

On the other hand, we have
$$
p_0\in \partial(\Int A)\subset \partial A = \partial (M\setminus A).
$$
As $M\setminus A$ is open in $M$, $p_0$ does not belong to $M\setminus A$ and $U\cap (M\setminus A)\neq \emptyset$.

Thus, any open subset containing $p_0$ includes an open subset where $\tau(\phi)$ vanishes everywhere and an open subset where $\tau(\phi)\neq 0$ at any point.

Let $V$ be an open subset containing $\phi(p_0)$ and $U$ an open subset containing $p_0$ such that $\phi(U)\subset V$.
Assume that $U$ and $V$ are the domains of local charts. Consider an open subset $D$ in $M$, containing $p_0$, such that its closure in $M$ is compact and contained in $U$.
As we have seen, the set $D$ contains an open subset where $\tau(\phi)$ vanishes everywhere and an open subset where $\tau(\phi)\neq 0$ at any point.

Denote
$$
\tau(\phi) = u^{\alpha}\frac{\partial}{\partial y^{\alpha}}.
$$
From \eqref{local-coordinates-lemma} we have that on $U$, and so on $D$,
\begin{align*}
\Delta u^{\theta} = & - 2 u^{\alpha}_j g^{ij} \phi^{\beta}_i \Gamma^{\theta}_{\beta \alpha} \\
& - u^{\alpha}\left( \left( \Delta \phi^{\beta} \right) \Gamma^{\theta}_{\beta \alpha} + g^{ij} \phi^{\beta}_j \phi^{\omega}_i \left(  \frac{\partial \Gamma^{\theta}_{\beta \omega}}{\partial y^{\alpha}} + \Gamma^{\gamma}_{\beta \omega}\Gamma^{\theta}_{\alpha \gamma}\right) \right).
\end{align*}
From the above equality we get the following inequality on $D$
\begin{align*}
\left|Au^{\alpha}\right|\leq C\left(\sum_{\beta,i}\left|\frac{\partial u^{\beta}}{\partial x^i}\right|+\sum_{\beta}\left|u^{\beta}\right|\right),
\end{align*}
as all functions $g^{ij}$, $\phi^{\alpha}$ and their derivatives, $\Gamma^{\theta}_{\beta \alpha}$ and their derivatives, are bounded on $D$.
Applying Theorem \ref{aro-theorem} we can deduce that \(u=0\), which implies that $\tau(\phi)$ vanishes everywhere in $D$. This contradiction implies that $\partial (\Int A)=\emptyset$ and we end the proof.
\end{proof}

\subsection{The case of a spherical target}

In this subsection we will prove Theorem \ref{sphere}.
First, we will exploit the fact that we are considering a spherical target
to bring \eqref{biharmonic} into a simpler form.

To this end, we consider the inclusion map \(\iota\colon \mathbb{S}^n\to\R^{n+1}\)
and form the composite map \(\varphi:=\iota\circ\phi\colon M\to\R^{n+1}\).

\begin{Lem}
For \(\varphi\colon M\to \mathbb{S}^n\subset\R^{n+1}\) with the constant curvature metric,
the equation for biharmonic maps acquires the form
\begin{align}
\label{biharmonic-spherical-target}
\Delta^2\varphi+&\left(|\Delta\varphi|^2+2|\nabla d\varphi|^2
+4\langle d\varphi,\nabla\Delta\varphi\rangle+2|d\varphi|^4+2\left\langle d\varphi\left(\operatorname{Ric}^M\right),d\varphi\right\rangle\right)\varphi\\
\nonumber&+4\sum_{i,j}^m\left\langle\left(\nabla d\varphi\right)\left(e_i,e_j\right),d\varphi\left(e_i\right)\right\rangle d\varphi\left(e_j\right)+2|d\varphi|^2\Delta\varphi=0,
\end{align}
where \(\operatorname{Ric}^M\) denotes the Ricci tensor field on the domain manifold \(M\)
and \(\{e_i\},i=1,\ldots,m\) is an orthonormal frame field.
\end{Lem}

\begin{proof}
It is well-known that for a spherical target of constant curvature \(1\) the following formula
holds true
\begin{align*}
d\iota(\tau(\phi))=\tau(\phi)=\Delta\varphi+|d\varphi|^2\varphi.
\end{align*}
Note that if \(\sigma\in\Gamma(\phi^\ast TN)\) then we can also think
of \(\sigma\) as a section in the pull-back bundle \(\varphi^{\ast}T\R^{n+1}\)
and the connections along \(\phi\) and \(\varphi\) are related via
\begin{align*}
\nabla_X^\varphi\sigma=\nabla_X^\phi\sigma-\langle d\phi(X),\sigma\rangle\varphi,\qquad\forall X\in\Gamma(TM).
\end{align*}

By a direct calculation one finds that
\begin{align*}
\tau_2(\phi)=&\Delta^2\varphi+2|d\varphi|^2\Delta\varphi+2d\varphi\left(\operatorname{grad}\left(|d\varphi|^2\right)\right)
+\left(\Delta|d\varphi|^2+2\operatorname{div}\theta^\sharp-|\Delta\varphi|^2+2|d\varphi|^4\right)\varphi,
\end{align*}
where \(\theta(X)=\langle d\phi(X),\tau(\phi)\rangle=\langle d\varphi(X),\Delta\varphi\rangle\),
see \cite{cezaryelin} for a derivation.

As a next step, we prove that
\begin{align*}
\operatorname{div}\theta^\sharp=|\Delta\varphi|^2+\langle d\varphi,\nabla\Delta\varphi\rangle,
\end{align*}
where \(\nabla=\nabla^\varphi\).
To this end we fix \(p\in M\) arbitrary and let \(\{X_i\},i=1,\ldots,m\) be a geodesic
frame field around \(p\). We calculate at the point \(p\)
\begin{align*}
\langle d\varphi,\nabla\Delta\varphi\rangle
=&\sum_{i=1}^m\left\langle d\varphi\left(X_i\right),\nabla_{X_i}\Delta\varphi\right\rangle\\
=&\sum_{i=1}^m\left(X_i\left\langle d\varphi\left(X_i\right),\Delta\varphi\right\rangle-\left\langle\nabla_{X_i}d\varphi\left(X_i\right),\Delta\varphi\right\rangle\right) \\
=&\sum_{i=1}^m\left(X_i\left\langle \theta^\sharp,X_i\right\rangle-\left\langle\left(\nabla d\varphi\right)\left(X_i,X_i\right),\Delta\varphi\right\rangle\right) \\
=&\sum_{i=1}^m\left(\left\langle\nabla_{X_i}\theta^\sharp,X_i\right\rangle-\left\langle\left(\nabla d\varphi\right)\left(X_i,X_i\right),\Delta\varphi\right\rangle\right) \\
=&\operatorname{div}\theta^\sharp-|\Delta\varphi|^2.
\end{align*}
To manipulate \(\Delta|d\varphi|^2\) we again compute at \(p\)
\begin{align*}
\Delta|d\varphi|^2=&\sum_{i=1}^mX_iX_i\langle d\varphi,d\varphi\rangle\\
=&2\sum_{i=1}^m\left(\left\langle\nabla_{X_i}\nabla_{X_i}d\varphi,d\varphi\right\rangle+\left|\nabla_{X_i}d\varphi\right|^2\right)\\
=&2\left\langle\tr_g\nabla^2d\varphi,d\varphi\right\rangle+2|\nabla d\varphi|^2.
\end{align*}
Recall that \(\tr_g\nabla^2d\varphi=-\Delta d\varphi+d\varphi(\operatorname{Ric}^M)\),
where \(\Delta\) is the Hodge-Laplacian acting on \(\R^{n+1}\)-valued one-forms.
Due to our sign convention we have \(-\Delta d\varphi=d\Delta\varphi=\nabla\Delta\varphi\)
such that we find
\begin{align*}
\Delta|d\varphi|^2=2|\nabla d\varphi|^2+2\langle\nabla\Delta\varphi,d\varphi\rangle+
2\left\langle d\varphi\left(\operatorname{Ric^M}\right),d\varphi\right\rangle.
\end{align*}
Finally, we have at the point \(p\)
\begin{align*}
d\varphi\left(\operatorname{grad}\left(|d\varphi|^2\right)\right)=2\sum_{i,j=1}^m\langle(\nabla d\varphi)(X_i,X_j),d\varphi(X_i)\rangle d\varphi(X_j)
\end{align*}
and the claim follows by combing the previous equations.
\end{proof}

Now, we define new variables \(v=\nabla\varphi\) and \(w=\Delta\varphi\).
In terms of a local chart \((U,x^i)\) on \(M\) and \(\{e_a\}\)
denoting the canonical basis of \(\R^{n+1}\), $a=1,\ldots,n+1$, we can write
\(v=\varphi_i^a dx^i\otimes e_a\), where \(\varphi_i^a=\frac{\partial\varphi^a}{\partial x^i}\).

Then any solution of \eqref{biharmonic-spherical-target} satisfies the second order elliptic
equation
\begin{align}
\label{def-F-sphere}
\Delta y=F(\varphi,\nabla\varphi,\nabla v,w,\nabla w),
\end{align}
where
\begin{align*}
y&=
\begin{pmatrix}
\varphi\\
v\\
w
\end{pmatrix},\\
F&=
\begin{pmatrix}
w \\
-dw \\
-\left( |w|^2+2|\nabla v|^2+4\langle v,\nabla w\rangle+2\sum_{i=1}^m \left\langle v \left( \operatorname{Ric}^M \left( e_i \right) \right),v \left( e_i \right) \right\rangle +2|v|^4 \right) \varphi \\
-4\sum_{i,j=1}^m \left\langle \left( \nabla v \right) \left( e_i,e_j \right),v\left( e_i \right) \right\rangle v \left( e_j\right)-2|v|^2w
\end{pmatrix}.
\end{align*}

At this point we consider two biharmonic maps \(\varphi_1,\varphi_2\) and their corresponding
new variables \(v_1,v_2,w_1,w_2\). We set \(u:=y_1-y_2\).
Note that the function \(u\) takes values in \(\R^{(n+1)(m+2)}\), which can be seen from
\begin{align*}
u&=
\begin{pmatrix}
\varphi_1-\varphi_2&=&\left( u^a \right)_a,\qquad a=1,\ldots,n+1 \\
v_1-v_2&=&\left( u^{n+1+\tilde a} \right)_{\tilde a},\qquad \tilde a=1,\ldots,m(n+1)\\
w_1-w_2&=&\left( u^{(n+1)(m+1)+a} \right)_a,\qquad a=1,\ldots,n+1.
\end{pmatrix}.
\end{align*}
Then we find
\begin{align*}
\Delta u&=
\begin{pmatrix}
w_1-w_2&=&\left( u^{(n+1)(m+1)+a}\right)_a\\
-\nabla \left( w_1-w_2 \right)&=&-\left( u_i^{(n+1)(m+1)+a} \right)_{i,a}\\
\Delta \left( w_1- w_2 \right)
\end{pmatrix},
\end{align*}
where we think of $u$ as being defined on the image of $U$ and $u_i^{(n+1)(m+1)+a}=\frac{\partial u^{(n+1)(m+1)+a}}{\partial x^i}$.

As a next step we prove that \(\Delta\left(w_1- w_2\right)\) can be expressed as a sum of
the components of \(u\) and their first order derivatives multiplied
by terms that do not contain the components of \(u\) or their derivatives. This allows us to obtain the estimate for $\left|\Delta\left(w_1-w_2\right)\right|$.

\begin{Lem}
Let \(D\subset U\) be an open subset such that its closure in $M$ is compact and included in $U$. Then the following estimate holds on $D$
\begin{align}
\label{estimate-sphere-w1w2}
\left|\Delta\left(w_1-w_2\right)\right|\leq C\left(\left|w_1-w_2\right|+\left|\nabla \left(w_1-w_2\right)\right|+\left|v_1-v_2\right|+\left|\nabla\left(v_1-v_2\right)\right|+\left|\varphi_1-\varphi_2\right|\right),
\end{align}
where the positive constant \(C\) depends on \(D\) and the derivatives of \(\varphi_1\) and \(\varphi_2\) up to third order.
\end{Lem}

\begin{proof}
Using \eqref{def-F-sphere} we find
\begin{align*}
\Delta (w_1-w_2)=&-\left(|w_1|^2+2|\nabla v_1|^2+4\langle v_1,\nabla w_1\rangle+2\sum_{i=1}^m\left\langle v_1\left(\operatorname{Ric}^M\left(e_i\right)\right),v_1\left(e_i\right)\right\rangle+2\left|v_1\right|^4\right)\varphi_1 \\
&-4\sum_{i,j=1}^m\langle(\nabla v_1)(e_i,e_j),v_1(e_i)\rangle v_1(e_j)-2|v_1|^2w_1 \\
&+\left(|w_2|^2+2|\nabla v_2|^2+4\langle v_2,\nabla w_2\rangle+2\sum_{i=1}^m\left\langle v_2\left(\operatorname{Ric}^M\left(e_i\right)\right),v_2\left(e_i\right)\right\rangle+2\left|v_2\right|^4\right)\varphi_2 \\
&+4\sum_{i,j=1}^m\langle(\nabla v_2)(e_i,e_j),v_2(e_i)\rangle v_2(e_j)+2|v_2|^2w_2.
\end{align*}
In the following we will rewrite all the terms on the right hand side by adding suitable zeros,
starting with
\begin{align*}
-|w_1|^2\varphi_1+|w_2|^2\varphi_2=-\langle w_1,w_1-w_2\rangle\varphi_1-\langle w_1-w_2,w_2\rangle\varphi_1-|w_2|^2(\varphi_1-\varphi_2).
\end{align*}
The next contribution can be manipulated as follows
\begin{align*}
-|\nabla v_1|^2\varphi_1+|\nabla v_2|^2\varphi_2
=&-\langle\nabla v_1,\nabla v_1-\nabla v_2\rangle\varphi_1-\langle\nabla v_1-\nabla v_2,\nabla v_2\rangle\varphi_1
-|\nabla v_2|^2(\varphi_1-\varphi_2).
\end{align*}
Moreover, we perform the following manipulation
\begin{align*}
-\langle v_1,\nabla w_1\rangle\varphi_1+\langle v_2,\nabla w_2\rangle\varphi_2
=&-\langle v_1-v_2,\nabla w_1\rangle\varphi_1-\langle v_2,\nabla w_1-\nabla w_2\rangle\varphi_1\\
&-\langle v_2,\nabla w_2\rangle(\varphi_1-\varphi_2).
\end{align*}
The next term can be rewritten as
\begin{align*}
-\sum_{i=1}^m&\left\langle v_1\left(\operatorname{Ric}^M\left(e_i\right)\right),v_1\left(e_i\right)\right\rangle
+\sum_{i=1}^m\left\langle v_2\left(\operatorname{Ric}^M\left(e_i\right)\right),v_2\left(e_i\right)\right\rangle \\
=&-\sum_{i=1}^m\left\langle \left(v_1-v_2\right)\left(\operatorname{Ric}^M\left(e_i\right)\right),v_1\left(e_i\right)\right\rangle\varphi_1 \\
&-\sum_{i=1}^m\left\langle v_2\left(\operatorname{Ric}^M\left(e_i\right)\right),\left(v_1-v_2\right)\left(e_i\right)\right\rangle\varphi_1 \\
&-\sum_{i=1}^m\left\langle v_2\left(\operatorname{Ric}^M\left(e_i\right)\right),v_2\left(e_i\right)\right\rangle\left(\varphi_1-\varphi_2\right).
\end{align*}
In addition, we have
\begin{align*}
-|v_1|^4\varphi_1+|v_2|^4\varphi_2=&-\left(|v_1|^2+|v_2|^2\right)\left(|v_1|^2-|v_2|^2\right)\varphi_1-|v_2|^4(\varphi_1-\varphi_2) \\
=&-\left(|v_1|^2+|v_2|^2\right)\langle v_1+v_2,v_1-v_2\rangle\varphi_1-|v_2|^4(\varphi_1-\varphi_2).
\end{align*}
Moreover, we find
\begin{align*}
- \sum_{i,j=1}^m&\langle(\nabla v_1)(e_i,e_j),v_1(e_i)\rangle v_1(e_j)+\sum_{i,j=1}^m(\langle\nabla v_2)(e_i,e_j),v_2(e_i)\rangle v_2(e_j)=\\
& - \sum_{i,j=1}^m\langle(\nabla v_1-\nabla v_2)(e_i,e_j),v_1(e_i)\rangle v_1(e_j)
- \sum_{i,j=1}^m\langle(\nabla v_2)(e_i,e_j),(v_1-v_2)(e_i)\rangle v_1(e_j)\\
& - \sum_{i,j=1}^m\langle(\nabla v_2)(e_i,e_j),v_2(e_i)\rangle(v_1-v_2)(e_j).
\end{align*}
Finally, we obtain
\begin{align*}
-|v_1|^2w_1+|v_2|^2w_2=-\langle v_1,v_1-v_2\rangle w_1-\langle v_1-v_2,v_2\rangle w_1-|v_2|^2(w_1-w_2).
\end{align*}
The result then follows directly since the closure of \(D\) is assumed to be compact and both \(\varphi_1,\varphi_2\)
are smooth by assumption.
\end{proof}

\begin{proof}[Proof of Theorem \ref{sphere}]
Making use of \eqref{estimate-sphere-w1w2} we get the following inequality
\begin{align*}
\left|Au^a\right|\leq C\left(\sum_{b,i}\left|\frac{\partial u^b}{\partial x^i}\right|+\sum_b\left|u^b\right|\right),
\end{align*}
where, here, $a,b=1,\ldots,(n+1)(m+2)$.

Applying Theorem \ref{aro-theorem} we can deduce that $u=0$ on $D$. We finish
the proof by denoting $A:=\{p\in M : \phi_1(p)=\phi_2(p)\}$ and using the same
arguments as in the proof of Theorem 1.2.
\end{proof}

\subsection{The case of a general target}

In this section we prove Theorem \ref{general-target}.

First, we consider a biharmonic map \(\phi\colon M\to N\)
and a local chart \((U,x^i)\) on \(M\) and \((V,y^\alpha)\) on \(N\) such that \(\phi(U)\subset V\).
In order to avoid any notational confusion the corresponding expression for \(\phi\)
in local coordinates will be denoted by \(\hat\phi\), that is
\begin{align*}
\hat\phi=\left( \phi^1,\ldots,\phi^n\right)=\left(\phi^\alpha\right)_\alpha.
\end{align*}

\begin{Lem}
Let \(\phi\colon M\to N\) be a biharmonic map. In terms of local coordinates it satisfies
\begin{align}
\label{biharmonic-general-target}
\Delta^2 \phi^{\theta} = & -4 \left\langle \nabla\Delta\phi^{\alpha},d\phi^{\beta} \right\rangle \Gamma^\theta_{\alpha\beta}
-2 \left\langle d\phi^{\alpha} \left( \operatorname{Ric}^M \right),d\phi^{\beta} \right\rangle \Gamma^\theta_{\alpha\beta}
-2 \left\langle \nabla d\phi^\alpha,\nabla d\phi^\beta \right\rangle \Gamma^\theta_{\alpha\beta} \\
\nonumber & -4g^{ij} \left\langle\nabla_{\frac{\partial}{\partial x^i}} d\phi^\alpha,d\phi^\beta \right\rangle \phi^\delta_jA^\theta_{\alpha\beta\delta}
 - 2\left\langle d\phi^\alpha,d\phi^\beta \right\rangle \left(\Delta\phi^\delta \right)A^\theta_{\alpha\beta\delta} \\
\nonumber &- \left\langle d\phi^\alpha,d\phi^\beta \right\rangle \left\langle d\phi^\delta,d\phi^\gamma \right\rangle B^\theta_{\alpha\beta\delta\gamma}
- \left(\Delta\phi^{\alpha}\right) \left(\Delta \phi^{\beta}\right) \Gamma^{\theta}_{\alpha \beta},
\end{align}
where
\begin{align*}
A^\theta_{\alpha\beta\delta}:=&\frac{\partial\Gamma^\theta_{\alpha\beta}}{\partial y^\delta}+\Gamma^\gamma_{\alpha\beta}\Gamma^\theta_{\gamma\delta},\\
B^\theta_{\alpha\beta\delta\gamma}:=&\frac{\partial^2\Gamma^\theta_{\alpha\beta}}{\partial y^\delta\partial y^\gamma}+\frac{\partial\Gamma^\omega_{\delta\gamma}}{\partial y^\alpha}\Gamma^\theta_{\omega\beta}
+\frac{\partial\Gamma^\omega_{\delta\gamma}}{\partial y^\beta}\Gamma^\theta_{\omega\alpha}
+\Gamma^\omega_{\delta\gamma}\frac{\partial\Gamma^\theta_{\alpha\beta}}{\partial y^\omega}+\Gamma^\omega_{\delta\gamma}\Gamma^\sigma_{\alpha\beta}\Gamma^\theta_{\omega\sigma}.
\end{align*}
\end{Lem}

\begin{proof}
Using the local form of the tension field \( u^{\theta}=\tau^\theta(\phi)=\Delta\phi^\theta+\left\langle d\phi^\alpha,d\phi^\beta\right\rangle\Gamma^\theta_{\alpha\beta} \)
in \eqref{local-coordinates-lemma} and the fact that
\begin{align*}
\tr\nabla^2 d\phi^\alpha=\nabla\Delta\phi^\alpha+d\phi^\alpha\left(\operatorname{Ric}^M\right),
\end{align*}
we get
\begin{align*}
0 = & \Delta^2 \phi^{\theta}+\Delta \left( \left\langle d\phi^\alpha,d\phi^\beta \right\rangle \Gamma^\theta_{\alpha\beta} \right) \\
& + 2 \left\langle d  \Delta\phi^{\alpha},d  \phi^{\beta} \right\rangle \Gamma^{\theta}_{\alpha \beta}
+ 2 \left\langle d  \left( \left\langle d\phi^\delta,d\phi^\gamma \right\rangle\Gamma^\alpha_{\delta\gamma} \right),d  \phi^{\beta}\right\rangle \Gamma^{\theta}_{\alpha \beta}  \\
& + \Delta\phi^{\alpha} \left( \left( \Delta \phi^{\beta} \right) \Gamma^{\theta}_{\alpha \beta}
+ \left\langle d \phi^{\beta},d \phi^{\omega} \right\rangle A^\theta_{\beta\omega\alpha} \right) \\
& + \left\langle d\phi^\delta,d\phi^\gamma \right\rangle \Gamma^\alpha_{\delta\gamma} \left( \left( \Delta \phi^{\beta} \right) \Gamma^{\theta}_{\alpha \beta}
+ \left\langle d \phi^{\beta},d \phi^{\omega} \right\rangle A^\theta_{\beta\omega\alpha} \right).
\end{align*}
By a direct calculation we obtain
\begin{align*}
\Delta\left( \left\langle d\phi^\alpha,d\phi^\beta \right\rangle\Gamma^\theta_{\alpha\beta} \right) = & 2 \left\langle \nabla\Delta\phi^{\alpha},d\phi^{\beta} \right\rangle \Gamma^\theta_{\alpha\beta}
+2 \left\langle d\phi^{\alpha} \left( \operatorname{Ric}^M \right),d\phi^{\beta} \right\rangle \Gamma^\theta_{\alpha\beta}
+2 \left\langle \nabla d\phi^\alpha,\nabla d\phi^\beta \right\rangle \Gamma^\theta_{\alpha\beta} \\
& +4g^{ij} \left\langle\nabla_{\frac{\partial}{\partial x^i}} d\phi^\alpha,d\phi^\beta \right\rangle \phi^\delta_j\frac{\partial\Gamma^\theta_{\alpha\beta}}{\partial y^\delta} \\
& + \left\langle d\phi^\alpha,d\phi^\beta \right\rangle \left(\Delta\phi^\delta \right)\frac{\partial\Gamma^\theta_{\alpha\beta}}{\partial y^\delta}
+ \left\langle d\phi^\alpha,d\phi^\beta \right\rangle \left\langle d\phi^\delta,d\phi^\gamma \right\rangle \frac{\partial^2\Gamma^\theta_{\alpha\beta}}{\partial y^\delta\partial y^\gamma}
\end{align*}
and also
\begin{align*}
\left\langle d  \left( \left\langle d\phi^\delta,d\phi^\gamma \right\rangle \Gamma^\alpha_{\delta\gamma} \right),d  \phi^{\beta} \right\rangle \Gamma^{\theta}_{\alpha \beta}
= &  2g^{ij} \left\langle \nabla_{\frac{\partial}{\partial x^i}} d\phi^\delta,d\phi^\gamma \right\rangle \phi^\beta_j \Gamma^\alpha_{\delta\gamma}\Gamma^\theta_{\alpha\beta} \\
& +\langle d\phi^\delta,d\phi^\gamma\rangle\langle d\phi^\omega,d\phi^\beta\rangle\frac{\partial\Gamma^\alpha_{\delta\gamma}}{\partial y^\omega}\Gamma^\theta_{\alpha\beta}.
\end{align*}
The claim follows by combining the equations.
\end{proof}

Now, as in the spherical case, we define new variables
\begin{align*}
v=\nabla\hat\phi=d\hat\phi=\left(d\phi^\alpha\right) \quad \textrm{and} \quad  w=\Delta\hat\phi=\left(\Delta\phi^\alpha\right).
\end{align*}
We can write
\(v=\phi_i^\alpha dx^i\otimes e_\alpha\), where \(\phi_i^\alpha=\frac{\partial\phi^\alpha}{\partial x^i}\) and \(\{e_\alpha\}\)
denotes the canonical basis of \(\R^{n}\).

Then any solution of \eqref{biharmonic-general-target} satisfies the second order elliptic
equation
\begin{align}
\label{def-F-general}
\Delta y=F(\hat\phi,\nabla\hat\phi,\nabla v,w,\nabla w),
\end{align}
where
\begin{align*}
y&=
\begin{pmatrix}
\hat\phi\\
v\\
w
\end{pmatrix},\qquad
F=
\begin{pmatrix}
w \\
-dw \\
F_3
\end{pmatrix}.
\end{align*}

Here, \(F_3\) is given by
\begin{align*}
F_3^\theta = & -4\left\langle \nabla w^\alpha,v^\beta \right\rangle \Gamma^\theta_{\alpha\beta}
-2 \left\langle v^\alpha\left(\operatorname{Ric}^M\right),v^\beta \right\rangle \Gamma^\theta_{\alpha\beta}
-2 \left\langle \nabla v^\alpha,\nabla v^\beta \right\rangle \Gamma^\theta_{\alpha\beta} \\
\nonumber & -4\left( \tr \left\langle \nabla_{(\cdot)} v^\alpha,v^\beta \right\rangle v^\delta(\cdot)\right)A^\theta_{\alpha\beta\delta}
 - 2\left\langle v^\alpha,v^\beta \right\rangle w^\delta A^\theta_{\alpha\beta\delta} \\
\nonumber &- \left\langle v^\alpha,v^\beta \right\rangle \left\langle v^\delta,v^\gamma \right\rangle B^\theta_{\alpha\beta\delta\gamma}
- w^{\alpha}w^{\beta} \Gamma^{\theta}_{\alpha \beta}.
\end{align*}

At this point we consider two biharmonic maps \(\phi_1,\phi_2\colon M\to N\) and
let \((U,x^i)\) be a local chart on \(M\) and \((V,y^\alpha)\) a local chart
on \(N\) such that \(\phi_1(U)\subset V\) and also \(\phi_2(U)\subset V\).

We set \(u:=y_1-y_2\).
Note that the function \(u\) takes values in \(\R^{n(m+2)}\), which can be seen from
\begin{align*}
u&=
\begin{pmatrix}
\hat\phi_1-\hat\phi_2&=&(u^\alpha)_\alpha,\qquad \alpha=1,\ldots,n \\
v_1-v_2&=&(u^{n+\tilde a})_{\tilde a},\qquad \tilde a=1,\ldots,mn\\
w_1-w_2&=&(u^{n(m+1)+\alpha})_\alpha,\qquad \alpha=1,\ldots,n
\end{pmatrix}.
\end{align*}
Then we find
\begin{align*}
\Delta u&=
\begin{pmatrix}
w_1-w_2&=&(u^{n(m+1)+\alpha})_\alpha\\
-\nabla(w_1-w_2)&=&-(u_i^{n(m+1)+\alpha})_{i,\alpha}\\
\Delta (w_1- w_2)
\end{pmatrix}.
\end{align*}

In order to obtain the estimate for $\left| \Delta \left( w_1 - w_2\right) \right|$, we need the image of $V$ in $\mathbb{R}^n$ to be convex and with compact closure
such that we can apply the mean-value theorem for functions defined on the image of $V$.
Then, when applying the mean value inequality, the standard norm of the differential of the function will be bounded on the image of $V$.
So, first, we consider the image of $V$ to be an open ball of radius $\varepsilon$ in the Euclidean space $\mathbb{R}^n$ and then we shrink it to a ball of radius $\frac{\varepsilon}{2}$.
Accordingly, we consider a smaller domain $U$. Thus, we will make use of the mean-value theorem for functions defined on the open ball of radius $\frac{\varepsilon}{2}$ but which are also defined on the closed ball of radius $\frac{\varepsilon}{2}$.

More precisely, for $f\colon \mathbb{B}^n\left(\frac{\varepsilon}{2}\right)\to \R$ we will apply the following inequality
\begin{align*}
|f(y_1) - f(y_2)|\leq |Df(y^\ast)||y_1 - y_2|,
\end{align*}
where \(y^\ast\) belongs to the standard segment that joins \(y_1\) with \(y_2\).

Then, as a next step, we prove that \(\Delta \left( w_1- w_2 \right) \) can be expressed as a sum of
the components of \(u\) and their first order derivatives multiplied
by terms that do not contain the components of \(u\) or their derivatives.

We have

\begin{Lem}
Let $D\subset U$ be an open subset such that its closure in $M$ is compact and included in $U$. Then the following estimate holds on $D$
\begin{align}
\label{estimate-w1w2-general}
\left|\Delta \left (w_1-w_2 \right) \right| \leq C \left( \left| w_1-w_2 \right| + \left| \nabla \left( w_1-w_2 \right) \right| + \left| v_1-v_2\right|+\left|\nabla\left( v_1-v_2 \right) \right|+\left| \hat\phi_1 - \hat\phi_2 \right| \right),
\end{align}
where the positive constant \(C\) depends on \(D\) and the derivatives of \(\hat\phi_1\) and \(\hat\phi_2\) up to third order.
\end{Lem}
\begin{proof}
Using \eqref{biharmonic-general-target} we find
\begin{align*}
\Delta w_1^\theta-\Delta w_2^\theta=&
- 4\left\langle \nabla w_1^\alpha,v_1^\beta \right\rangle \Gamma(\hat\phi_1)^\theta_{\alpha\beta}
+ 4\left\langle \nabla w_2^\alpha,v_2^\beta \right\rangle \Gamma(\hat\phi_2)^\theta_{\alpha\beta}
\\
\nonumber & -2 \left\langle v_1^\alpha\left(\operatorname{Ric}^M\right),v_1^\beta \right\rangle \Gamma(\hat\phi_1)^\theta_{\alpha\beta}
+ 2 \left\langle v_2^\alpha\left(\operatorname{Ric}^M\right),v_2^\beta \right\rangle \Gamma(\hat\phi_2)^\theta_{\alpha\beta}
\\
\nonumber & -2 \left\langle \nabla v_1^\alpha,\nabla v_1^\beta \right\rangle \Gamma(\hat\phi_1)^\theta_{\alpha\beta}
+2 \left\langle \nabla v_2^\alpha,\nabla v_2^\beta \right\rangle \Gamma(\hat\phi_2)^\theta_{\alpha\beta}
\\
\nonumber & - 4 \left( \tr \left\langle \nabla_{(\cdot)} v_1^\alpha,v_1^\beta \right\rangle v_1^\delta(\cdot)\right) A(\hat\phi_1)^\theta_{\alpha\beta\delta}
+ 4 \left( \tr \left\langle \nabla_{(\cdot)} v_2^\alpha,v_2^\beta \right\rangle v_2^\delta(\cdot)\right) A(\hat\phi_2)^\theta_{\alpha\beta\delta}
\\
\nonumber & - 2 \left\langle v_1^\alpha,v_1^\beta \right\rangle w_1^\delta A(\hat\phi_1)^\theta_{\alpha\beta\delta}
+ 2 \left\langle v_2^\alpha,v_2^\beta \right\rangle w_2^\delta A(\hat\phi_2)^\theta_{\alpha\beta\delta}
\\
\nonumber & - \left\langle v_1^\alpha,v_1^\beta \right\rangle \left\langle v_1^\delta,v_1^\gamma \right\rangle B(\hat\phi_1)^\theta_{\alpha\beta\delta\gamma}
+ \left\langle v_2^\alpha,v_2^\beta \right\rangle \left\langle v_2^\delta,v_2^\gamma \right\rangle B(\hat\phi_2)^\theta_{\alpha\beta\delta\gamma}
\\
\nonumber & - w_1^\alpha w_1^\beta \Gamma(\hat\phi_1)^{\theta}_{\alpha \beta}
+ w_2^\alpha w_2^\beta \Gamma(\hat\phi_2)^{\theta}_{\alpha \beta}.
\end{align*}
Now, we have to start estimating all the terms on the right hand side by adding suitable zeros.

The first term can be controlled as follows
\begin{align*}
&-\left\langle \nabla w_1^\alpha,v_1^\beta \right\rangle \Gamma(\hat\phi_1)^\theta_{\alpha\beta}
+\left\langle \nabla w_2^\alpha,v_2^\beta \right\rangle \Gamma(\hat\phi_2)^\theta_{\alpha\beta} \\
=&-\left(\left\langle \nabla w_1^\alpha-\nabla w_2^\alpha,v_1^\beta \right\rangle+
\left\langle\nabla w_2^\alpha,v_1^\beta-v_2^\beta\right\rangle\right)\Gamma(\hat\phi_1)^\theta_{\alpha\beta} \\
&-\left\langle \nabla w_2^\alpha,v_2^\beta \right\rangle\left(\Gamma(\hat\phi_1)^\theta_{\alpha\beta}-\Gamma(\hat\phi_2)^\theta_{\alpha\beta}\right)\\
\leq &C\left(\left|\nabla\left(w_1-w_2\right)\right|+\left|v_1-v_2\right|+\left|\hat\phi_1-\hat\phi_2\right|\right).
\end{align*}

For the second term we rewrite
\begin{align*}
&-\left\langle v_1^\alpha\left(\operatorname{Ric}^M\right),v_1^\beta\right\rangle\Gamma(\hat\phi_1)^\theta_{\alpha\beta}
+\left\langle v_2^\alpha\left(\operatorname{Ric}^M\right),v_2^\beta\right\rangle\Gamma(\hat\phi_2)^\theta_{\alpha\beta} \\
=&-\left(\left\langle \left(v_1-v_2\right)^\alpha\left(\operatorname{Ric}^M\right),v_1^\beta\right\rangle
+\left\langle v_2^\alpha\left(\operatorname{Ric}^M\right),\left(v_1-v_2\right)^\beta\right\rangle\right)\Gamma(\hat\phi_1)^\theta_{\alpha\beta} \\
&-\left\langle v_2^\alpha\left(\operatorname{Ric}^M\right),v_2^\beta\right\rangle\left(\Gamma(\hat\phi_1)^\theta_{\alpha\beta}-\Gamma(\hat\phi_2)^\theta_{\alpha\beta}\right) \\
\leq& C\left(\left|v_1-v_2\right|+\left|\hat\phi_1-\hat\phi_2\right|\right),
\end{align*}
where we applied the mean value inequality to the last term.
In the following we will frequently apply the mean value inequality without mentioning it explicitly.

The third contribution can be manipulated as follows
\begin{align*}
& -\left\langle \nabla v_1^\alpha,\nabla v_1^\beta \right\rangle \Gamma(\hat\phi_1)^\theta_{\alpha\beta}
+\left\langle \nabla v_2^\alpha,\nabla v_2^\beta \right\rangle \Gamma(\hat\phi_2)^\theta_{\alpha\beta} \\
=&-\left(\left\langle \nabla v_1^\alpha-\nabla v_2^\alpha,\nabla v_1^\beta \right\rangle
+\left\langle \nabla v_2^\alpha,\nabla v_1^\beta-\nabla v_2^\beta \right\rangle\right)\Gamma(\hat\phi_1)^\theta_{\alpha\beta}
-\left\langle \nabla v_2^\alpha,\nabla v_2^\beta\right\rangle\left(\Gamma(\hat\phi_1)^\theta_{\alpha\beta}-\Gamma(\hat\phi_2)^\theta_{\alpha\beta}\right)\\
\leq &C\left(\left|\nabla(v_1-v_2)\right|+\left|\hat\phi_1-\hat\phi_2\right|\right).
\end{align*}

Regarding the fourth term we find
\begin{align*}
& - \left( \tr \left\langle \nabla_{(\cdot)} v_1^\alpha,v_1^\beta \right\rangle v_1^\delta(\cdot)\right) A(\hat\phi_1)^\theta_{\alpha\beta\delta}
+ \left( \tr \left\langle \nabla_{(\cdot)} v_2^\alpha,v_2^\beta \right\rangle v_2^\delta(\cdot)\right) A(\hat\phi_2)^\theta_{\alpha\beta\delta}\\
=&-\left(
\left(\tr \left\langle \nabla_{(\cdot)} v_1^\alpha-\nabla_{(\cdot)}v_2^\alpha,v_1^\beta \right\rangle v_1^\delta(\cdot)\right)
+\left(\tr \left\langle \nabla_{(\cdot)} v_2^\alpha,v_1^\beta-v_2^\beta \right\rangle v_1^\delta(\cdot)\right) \right. \\
&\left. +\left(\tr \left\langle \nabla_{(\cdot)} v_2^\alpha,v_2^\beta\right\rangle \left(v_1^\delta(\cdot)-v_2^\delta(\cdot)\right)\right)
\right) A(\hat\phi_1)^\theta_{\alpha\beta\delta}\\
&-\left( \tr \left\langle \nabla_{(\cdot)} v_2^\alpha,v_2^\beta \right\rangle v_2^\delta(\cdot)\right)
\left(A(\hat\phi_1)^\theta_{\alpha\beta\delta}-A(\hat\phi_2)^\theta_{\alpha\beta\delta}\right) \\
\leq & C\left(\left|v_1-v_2\right|+\left|\nabla(v_1-v_2)\right|+\left|\hat\phi_1-\hat\phi_2\right|\right).
\end{align*}

For the fifth term we obtain
\begin{align*}
&- \left\langle v_1^\alpha,v_1^\beta \right\rangle w_1^\delta A(\hat\phi_1)^\theta_{\alpha\beta\delta}
+  \left\langle v_2^\alpha,v_2^\beta \right\rangle w_2^\delta A(\hat\phi_2)^\theta_{\alpha\beta\delta}\\
=&-\left(\left\langle (v_1-v_2)^\alpha,v_1^\beta\right\rangle w_1^\delta
+\left\langle v_2^\alpha,(v_1-v_2)^\beta\right\rangle w_1^\delta
+\left\langle v_2^\alpha,v_2^\beta\right\rangle (w_1-w_2)^\delta\right) A(\hat\phi_1)^\theta_{\alpha\beta\delta} \\
&-\left\langle v_2^\alpha,v_2^\beta\right\rangle w_2^\delta
\left(A(\hat\phi_1)^\theta_{\alpha\beta\delta}-A(\hat\phi_2)^\theta_{\alpha\beta\delta}\right) \\
\leq& C\left(\left|v_1-v_2\right|+\left|w_1-w_2\right|+\left|\hat\phi_1-\hat\phi_2\right|\right).
\end{align*}

The sixth contribution can be estimated as
\begin{align*}
\nonumber & - \left\langle v_1^\alpha,v_1^\beta \right\rangle \left\langle v_1^\delta,v_1^\gamma \right\rangle B(\hat\phi_1)^\theta_{\alpha\beta\delta\gamma}
+ \left\langle v_2^\alpha,v_2^\beta \right\rangle \left\langle v_2^\delta,v_2^\gamma \right\rangle B(\hat\phi_2)^\theta_{\alpha\beta\delta\gamma} \\
=&-\left(\left\langle (v_1-v_2)^\alpha,v_1^\beta\right\rangle\left\langle v_1^\delta,v_1^\gamma\right\rangle
+\left\langle v_2^\alpha,(v_1-v_2)^\beta\right\rangle\left\langle v_1^\delta,v_1^\gamma\right\rangle \right. \\
&\left. +\left\langle v_2^\alpha,v_2^\beta\right\rangle\left\langle (v_1-v_2)^\delta,v_1^\gamma\right\rangle
+\left\langle v_2^\alpha,v_2^\beta\right\rangle \left\langle v_2^\delta,(v_1-v_2)^\gamma\right\rangle\right)B(\hat\phi_1)^\theta_{\alpha\beta\delta\gamma} \\
&-\left\langle v_2^\alpha,v_2^\beta\right\rangle\left\langle v_2^\delta,v_2^\gamma\right\rangle
\left(B(\hat\phi_1)^\theta_{\alpha\beta\delta\gamma} -B(\hat\phi_2)^\theta_{\alpha\beta\delta\gamma} \right)\\
\leq& C\left(\left|v_1-v_2\right|+\left|\hat\phi_1-\hat\phi_2\right|\right).
\end{align*}

The seventh term can be estimated as
\begin{align*}
-w_1^\alpha w_1^\beta& \Gamma(\hat\phi_1)^{\theta}_{\alpha \beta}+w_2^\alpha w_2^\beta \Gamma(\hat\phi_2)^{\theta}_{\alpha \beta} \\
=&-\left((w_1-w_2)^\alpha w_1^\beta
+w_2^\alpha (w_1-w_2)^\beta \right)\Gamma(\hat\phi_1)^{\theta}_{\alpha \beta}
-w_2^\alpha w_2^\beta \left(\Gamma(\hat\phi_1)^{\theta}_{\alpha \beta}-\Gamma(\hat\phi_2)^{\theta}_{\alpha \beta}\right)\\
\leq &C\left(\left|w_1-w_2\right|+\left|\hat\phi_1-\hat\phi_2\right|\right).
\end{align*}
The claim follows by adding up the different contributions.
\end{proof}

\begin{proof}[Proof of Theorem \ref{sphere}]
Making use of \eqref{estimate-w1w2-general} we get the following inequality
\begin{align*}
\left|Au^a\right|\leq C\left(\sum_{b,i}\left|\frac{\partial u^b}{\partial x^i}\right|+\sum_b\left|u^b\right|\right),
\end{align*}
where $a,b=1,\ldots, n(m+2)$.

Applying Theorem \ref{aro-theorem} we can deduce that $u=0$ on $D$. We finish
the proof by denoting $A:=\{p\in M : \phi_1(p)=\phi_2(p)\}$ and using the same
arguments as in the proof of Theorem \ref{harmonic-everywhere}.
\end{proof}

\subsection{Proof of Theorem \ref{theorem-totally-geodesic-sphere}}

In this section we prove Theorem \ref{theorem-totally-geodesic-sphere}.
The proof is again based on Theorem \ref{aro-theorem} and the concrete
expressions of the Christoffel symbols on \(\mathbb{S}^n\).

Let \(\mathbb{S}^n\) be the Euclidean unit sphere and denote by \(N\) and \(S\)
the north and south pole, respectively.
It is well-known that
\begin{align*}
\mathbb{S}^n\setminus\{N,S\}=\left(\mathbb{S}^{n-1}\times (0,\pi),\sin^2 s \cdot g_{\mathbb{S}^{n-1}}+ds^2\right).
\end{align*}

Let \((y^a)\) be local coordinates on \(\mathbb{S}^{n-1},a=1,\ldots,n-1\).
Then \((y^1,\ldots,y^{n-1},y^n=s)\) are local coordinates on \(\mathbb{S}^n\setminus\{N,S\}\).

In this geometric setup the Christoffel symbols on \(\mathbb{S}^n\) are given by
\begin{align*}
\Gamma^a_{bc}(y^\alpha)=&\tilde\Gamma^a_{bc}(y^d),\\
\Gamma^n_{bc}(y^\alpha)=&-\sin s\cos s \ \tilde g_{bc}(y^d)=-\frac{1}{2}\sin(2s) \ \tilde g_{bc}(y^d),\\
\Gamma^a_{nn}(y^\alpha)=&\Gamma^n_{nn}(y^\alpha)=0,\\
\Gamma^a_{bn}(y^\alpha)=&\frac{\cos s}{\sin s}\delta^a_b,\\
\Gamma^n_{bn}(y^\alpha)=&0,
\end{align*}
where we use a \(\tilde{}\) to indicate objects on \(\mathbb{S}^{n-1}\).

The equator \(\mathbb{S}^{n-1}\) is given by
\begin{align*}
\mathbb{S}^{n-1}: y^n=s=\frac{\pi}{2}.
\end{align*}

Now, let \((U,x^i)\) be a local chart on \(M\) and we denote the domain of the above
local coordinates on \(\mathbb{S}^n\setminus\{N,S\}\) by \(V\).
In addition, we assume that \(\phi(U)\subset V\).

We denote the corresponding expression for \(\phi\) in local coordinates by \(\hat\phi\), i.e.
\begin{align*}
\hat\phi=(\phi^1,\ldots,\phi^n)=(\phi^\alpha).
\end{align*}

Assume that \(W\) is an open subset of \(U\) and \(\phi(W)\subset \mathbb{S}^{n-1}\), i.e. \(\phi(W)\subset \mathbb{S}^{n-1}\cap V\).
Hence, in \(W\) we have \(\phi^n=\frac{\pi}{2}\). Now, define \(f\colon U\to\R, f:=\phi^n-\frac{\pi}{2}\).
Clearly, \(f\) vanishes when restricted to \(W\).

Let \(D\) be an open subset of \(U\) such that its closure in \(M\)
is compact and included in $U$, and \(W\subset D\subset U\).

We will prove the following estimate:
\begin{Lem}
The function \(f\colon U\to\R\) defined above satisfies the following estimate on $D$
\begin{align}
|\Delta^2 f|\leq C\left(\left|\nabla\Delta f\right|+\left|\Delta f\right|+\left|\nabla df\right|+\left|\nabla f\right|+\left|f\right|\right).
\end{align}
\end{Lem}
\begin{proof}
We have $\Delta^2f=\Delta^2\phi^n$ and we will use ~\eqref{biharmonic-general-target} with $\theta=n$ in order to obtain the estimates.

We expand the first term as
\begin{align*}
\left\langle \nabla\Delta \phi^{\alpha},d\phi^{\beta}\right\rangle\Gamma^n_{\alpha\beta} = &
\left\langle \nabla\Delta \phi^{n},d\phi^{\beta}\right\rangle\Gamma^n_{n\beta} +
\left\langle \nabla\Delta \phi^{b},d\phi^{c}\right\rangle\Gamma^n_{bc} +
\left\langle \nabla\Delta \phi^{b},d\phi^{n}\right\rangle\Gamma^n_{bn} \\
\nonumber = & \left\langle \nabla\Delta f,d\phi^{\beta}\right\rangle\Gamma^n_{n\beta} +
\left\langle \nabla\Delta \phi^{b},d\phi^{c}\right\rangle\Gamma^n_{bc} +
\left\langle \nabla\Delta \phi^{b},\nabla f\right\rangle\Gamma^n_{bn}. \\
\end{align*}
Since
\begin{align*}
\left| \left\langle \nabla\Delta f,d\phi^{\beta}\right\rangle\Gamma^n_{n\beta} +
\left\langle \nabla\Delta \phi^{b},\nabla f\right\rangle\Gamma^n_{bn} \right| \leq C\left(\left|\nabla\Delta f\right| + \left|\nabla f\right|\right),
\end{align*}
it remains to estimate the term $\left| \left\langle \nabla\Delta \phi^{b},d\phi^{c}\right\rangle\Gamma^n_{bc} \right|$.
We immediately obtain
\begin{align*}
\left|\Gamma^n_{bc}\right| = & \left|-\frac{1}{2}\sin(2\phi^n)\tilde{g}_{bc} \right| = \left|\frac{1}{2}\sin(2f)\tilde{g}_{bc} \right| \\
\nonumber \leq & C|f|.
\end{align*}
Therefore, for the first term on the right hand side of \eqref{biharmonic-general-target} we have the estimate
\begin{align*}
\left|\left\langle \nabla\Delta \phi^{\alpha},d\phi^{\beta}\right\rangle\Gamma^n_{\alpha\beta}\right| \leq C \left( \left|\nabla\Delta f\right| + \left| \nabla f\right| + \left| f\right|\right).
\end{align*}
The second, third and last term can be estimated in the same manner as the first one.

We write the fourth term in the same way as the first one and, the only term we must estimate is
\begin{align*}
\left| g^{ij}\left\langle\nabla_{\frac{\partial}{\partial x^i}}d\phi^a,d\phi^b\right\rangle\phi^c_j A^n_{abc} \right|.
\end{align*}
We directly find
\begin{align*}
A^n_{abc}=\frac{1}{2}\sin(2f)\left( \frac{\partial \tilde{g}_{ab}}{\partial y^c} + \tilde{\Gamma}^e_{ab}\tilde{g}_{ec} \right),
\end{align*}
such that
\begin{align*}
\left| A^n_{abc}\right| \leq C|f|
\end{align*}
and
\begin{align*}
\left| g^{ij}\left\langle\nabla_{\frac{\partial}{\partial x^i}}d\phi^\alpha,d\phi^\beta\right\rangle\phi^\delta_j A^n_{\alpha\beta\delta} \right|
\leq C \left( \left|\nabla df\right| + \left|\nabla f\right| + \left| f\right| \right).
\end{align*}
The fifth and sixth terms will be estimated in the same way as the fourth term,
which completes the proof.
\end{proof}

\begin{proof}[Proof of Theorem \ref{theorem-totally-geodesic-sphere}]
As in the proofs of Theorems \ref{sphere} and \ref{general-target} we
define suitable functions \(u\) and $F$. Applying Theorem \ref{aro-theorem} we can deduce that $u=0$ on $D$, i.e. $\phi$ maps the whole of $D$ into $\mathbb{S}^{n-1}$.

We finish the proof by denoting $A:=\{p\in M : \phi(p)\in \mathbb{S}^{n-1}\}$ and using the same
arguments as in the proof of Theorem \ref{sampson-submanifold}.
\end{proof}

\subsection{Proof of Theorem \ref{theorem-totally-geodesic}}
In this section we prove Theorem \ref{theorem-totally-geodesic}.
Let \((V,y^\alpha)\) be a local chart on \(N\) such that
\begin{align*}
V\cap P\colon y^{r+1}=\ldots=y^n=0,
\end{align*}
where \(r\) denotes the dimension of the submanifold \(P\).
Moreover, we assume that the image of \(V\) in \(\R^n\) is convex.
For example, we can assume that it is the interior of a \(n\)-cube with side lengths \(\epsilon/2\)
that lies inside a \(n\)-cube with side lengths \(\epsilon\).

In addition, let \((U,x^i)\) be a local chart on \(M\) and assume that \(\phi(U)\subset V\).
We denote the corresponding expression for \(\phi\) in local coordinates by \(\hat\phi\), i.e.
\begin{align*}
\hat\phi=(\phi^1,\ldots,\phi^n)=(\phi^\alpha).
\end{align*}
As \(P\) is totally geodesic in \(N\), we have on the intersection \(V\cap P\)
\begin{align}
\label{christoffel-a}
\begin{cases}
\Gamma^\theta_{\alpha\beta}={}^P\Gamma^\theta_{\alpha\beta}&  \text{for all}~~ 1\leq\alpha,\beta\leq r, 1\leq\theta\leq r,\\
\Gamma^\theta_{\alpha\beta}=0&  \text{for all}~~ 1\leq\alpha,\beta\leq r, r+1\leq\theta\leq n,
\end{cases}
\end{align}
where \({}^P\Gamma^\theta_{\alpha\beta}\) denote the Christoffel symbols of the submanifold \(P\)
for \(1\leq\theta,\alpha,\beta\leq r\). On \(V\cap P\) we also have
\begin{align}
\label{christoffel-b}
\frac{\partial\Gamma^\theta_{\alpha\beta}}{\partial y^\delta}=\frac{\partial^2\Gamma^\theta_{\alpha\beta}}{\partial y^\delta\partial y^\sigma}=0,\qquad
\text{for all}~~ 1\leq\alpha,\beta,\delta,\sigma\leq r,~~ r+1\leq\theta\leq n.
\end{align}
Now, assume that \(W\) is an open subset of \(U\) and \(\phi(W)\subset P\).
Hence, on \(W\) we have
\begin{align*}
\phi^{r+1}=\ldots=\phi^n=0.
\end{align*}
We set
\begin{align*}
f=(\phi^{r+1},\ldots,\phi^n)=(f^1,\ldots,f^{n-r})=(f^{\tilde{a}}),
\end{align*}
where \(\tilde{a}=1,\ldots,n-r\). Clearly, we have \(f\colon U\to\R^{n-r}\) and also \(f|_{W}=0\).
Now, let \(D\) be an open subset of \(U\) such that its closure in \(M\) is compact and included in \(U\), and \(W\subset D\subset U\).

We will prove the following estimate:
\begin{Lem}
Restricted to \(D\), the function \(f\colon U\to\R^{n-r}\) defined above satisfies the following estimate
\begin{align*}
|\Delta^2f^{\tilde{a}}|\leq C(|\nabla\Delta f|+|\Delta f|+|\nabla df|+|\nabla f|+|f|),
\end{align*}
where the constant \(C\) depends on the geometry of \(N\).
\end{Lem}
\begin{proof}
We make use of \eqref{biharmonic-general-target} and estimate each term on the right-hand side.
We expand the first term of \eqref{biharmonic-general-target} in the following way
\begin{align*}
\langle\nabla\Delta\phi^\alpha,d\phi^\beta\rangle\Gamma_{\alpha\beta}^{r+\tilde{a}}
=&\langle\nabla\Delta\phi^{r+\tilde{b}},d\phi^\beta\rangle\Gamma_{({r+\tilde{b}})\beta}^{r+\tilde{a}}
+\langle\nabla\Delta\phi^b,d\phi^c\rangle\Gamma_{bc}^{r+\tilde{a}}
+\langle\nabla\Delta\phi^b,d\phi^{r+\tilde{c}}\rangle\Gamma_{b(r+\tilde{c})}^{r+\tilde{a}}\\
=&\langle\nabla\Delta f^{\tilde{b}},d\phi^\beta\rangle\Gamma_{({r+\tilde{b}})\beta}^{r+\tilde{a}}
+\langle\nabla\Delta\phi^b,d\phi^c\rangle\Gamma_{bc}^{r+\tilde{a}}
+\langle\nabla\Delta\phi^b,d f^{\tilde{c}}\rangle\Gamma_{b(r+\tilde{c})}^{r+\tilde{a}},
\end{align*}
where \(b,c=1,\ldots,r\).
Since we have
\begin{align*}
\left|\langle\nabla\Delta f^{\tilde{b}},d\phi^\beta\rangle\Gamma_{({r+\tilde{b}})\beta}^{r+\tilde{a}}
+\langle\nabla\Delta\phi^b,d f^{\tilde{c}}\rangle\Gamma_{b(r+\tilde{c})}^{r+\tilde{a}}\right|
\leq C(|\nabla\Delta f|+|\nabla f|),
\end{align*}
it remains to estimate the term \(\langle\nabla\Delta\phi^b,d\phi^c\rangle\Gamma_{bc}^{r+\tilde{a}}\).
To this end we define a new function
\(\hat\phi'\colon U\to\R^n\) by
\begin{align*}
\hat\phi'=(\phi^1,\ldots,\phi^r,0,\ldots,0)
\end{align*}
and consequently \(\phi'\colon U\to N\) takes values in \(P\).
We note that, according to \eqref{christoffel-a}, we have \(\Gamma_{bc}^{r+\tilde{a}}=\Gamma_{bc}^{r+\tilde{a}}(\phi)=0\)
on \(W\), but at the points in \(D\setminus W\) the image of the map \(\phi\) may not be in \(P\)
such that \(\Gamma_{bc}^{r+\tilde{a}}\neq 0\).
However, on \(U\), and thus also on \(D\), we have
\begin{align*}
\Gamma_{bc}^{r+\tilde{a}}=&\Gamma_{bc}^{r+\tilde{a}}(\phi)-0 \\
=&\Gamma_{bc}^{r+\tilde{a}}(\phi)-\Gamma_{bc}^{r+\tilde{a}}(\phi').
\end{align*}
By applying the mean value inequality we get
\begin{align}
\label{christoffel-c}
\left|\Gamma_{bc}^{r+\tilde{a}}\right|\leq C\left|f\right|.
\end{align}
Consequently, for the first term on the right hand side of \eqref{biharmonic-general-target},
we get the estimate
\begin{align*}
\left|\langle\nabla\Delta\phi^\alpha,d\phi^\beta\rangle\Gamma_{\alpha\beta}^{r+\tilde{a}}\right|\leq C(|\nabla\Delta f|+|\nabla f|+|f|).
\end{align*}
The second, third and the last term of \eqref{biharmonic-general-target} can be estimated in the same manner as the first one.
In order to estimate the fourth term we rewrite it in the same way as the first one, and, the only
contribution we have to estimate is
\begin{align*}
\left|g^{ij} \left\langle\nabla_{\frac{\partial}{\partial x^i}} d\phi^a,d\phi^b \right\rangle \phi^c_jA^{r+\tilde{a}}_{abc}\right|.
\end{align*}
Remember that
\begin{align*}
A^{r+\tilde{a}}_{abc}=\frac{\partial\Gamma^{r+\tilde{a}}_{ab}}{\partial y^c}+\Gamma^\gamma_{ab}\Gamma^{r+\tilde{a}}_{\gamma c}.
\end{align*}
According to \eqref{christoffel-b} and \eqref{christoffel-c} we get on \(D\)
\begin{align}
\label{christoffel-d}
\left|\frac{\partial\Gamma^{r+\tilde{a}}_{ab}}{\partial y^c}\right|\leq C|f|.
\end{align}
Moreover, taking into account \eqref{christoffel-c}, we obtain
\begin{align*}
\left|\Gamma^\gamma_{ab}\Gamma^{r+\tilde{a}}_{\gamma c}\right|
=\left|\Gamma^d_{ab}\Gamma^{r+\tilde{a}}_{dc}+|\Gamma^{r+\tilde{d}}_{ab}\Gamma^{r+\tilde{a}}_{(r+\tilde{d})c}
\right|\leq C|f|.
\end{align*}
Thus, the fourth term of \eqref{biharmonic-general-target} can also be estimated such that it gives the desired inequality.
The fifth term of \eqref{biharmonic-general-target} can be estimated in the same way as the fourth one.
Regarding the sixth term we have to estimate
\begin{align*}
B^{r+\tilde{a}}_{abcd}=&\frac{\partial^2\Gamma^{r+\tilde{a}}_{ab}}{\partial y^c\partial y^d}
+\frac{\partial\Gamma^\omega_{cd}}{\partial y^a}\Gamma^{r+\tilde{a}}_{\omega b}
+\frac{\partial\Gamma^\omega_{cd}}{\partial y^b}\Gamma^{r+\tilde{a}}_{\omega a}
+\Gamma^\omega_{cd}\frac{\partial\Gamma^{r+\tilde{a}}_{ab}}{\partial y^\omega}
+\Gamma^\omega_{cd}\Gamma^\sigma_{ab}\Gamma^{r+\tilde{a}}_{\omega\sigma}.
\end{align*}
Applying \eqref{christoffel-b} and similar as in \eqref{christoffel-c} we get
\begin{align*}
\left|\frac{\partial^2\Gamma^{r+\tilde{a}}_{ab}}{\partial y^c\partial y^d}\right|\leq C|f|.
\end{align*}
Then, making use of \eqref{christoffel-c} and \eqref{christoffel-d} we obtain
\begin{align*}
\left|\frac{\partial\Gamma^\omega_{cd}}{\partial y^a}\Gamma^{r+\tilde{a}}_{\omega b}\right|
=\left|\frac{\partial\Gamma^l_{cd}}{\partial y^a}\Gamma^{r+\tilde{a}}_{lb}
+\frac{\partial\Gamma^{r+\tilde{l}}_{cd}}{\partial y^a}\Gamma^{r+\tilde{a}}_{(r+\tilde{l}) b}
\right|
\leq C|f|.
\end{align*}
In the same manner we find
\begin{align*}
\left|\frac{\partial\Gamma^\omega_{cd}}{\partial y^b}\Gamma^{r+\tilde{a}}_{\omega a}\right|
\leq C|f|.
\end{align*}
In addition, we get
\begin{align*}
\left|\Gamma^\omega_{cd}\frac{\partial\Gamma^{r+\tilde{a}}_{ab}}{\partial y^\omega}\right|
=\left|\Gamma^l_{cd}\frac{\partial\Gamma^{r+\tilde{a}}_{ab}}{\partial y^l}
+\Gamma^{r+\tilde{l}}_{cd}\frac{\partial\Gamma^{r+\tilde{a}}_{ab}}{\partial y^{r+\tilde{l}}}
\right|
\leq C|f|.
\end{align*}
Finally, we have
\begin{align*}
\left|\Gamma^\omega_{cd}\Gamma^\sigma_{ab}\Gamma^{r+\tilde{a}}_{\omega\sigma}\right|
=&\left|\Gamma^l_{cd}\Gamma^\sigma_{ab}\Gamma^{r+\tilde{a}}_{l\sigma}
+\Gamma^{r+\tilde{l}}_{cd}\Gamma^\sigma_{ab}\Gamma^{r+\tilde{a}}_{(r+\tilde{l})\sigma}
\right|\\
=&\left|\Gamma^l_{cd}\Gamma^e_{ab}\Gamma^{r+\tilde{a}}_{le}
+\Gamma^l_{cd}\Gamma^{r+\tilde{e}}_{ab}\Gamma^{r+\tilde{a}}_{l(r+\tilde{e})}
+\Gamma^{r+\tilde{l}}_{cd}\Gamma^\sigma_{ab}\Gamma^{r+\tilde{a}}_{(r+\tilde{l})\sigma}
\right|
\leq C|f|,
\end{align*}
which establishes the desired estimate.
\end{proof}

\begin{proof}[Proof of Theorem \ref{theorem-totally-geodesic}]
As in the proofs of Theorems \ref{sphere} and \ref{general-target} we
define suitable functions \(u\) and $F$. Applying Theorem \ref{aro-theorem} we can deduce that $u=0$ on $D$, i.e. $\phi$ maps the whole of $D$ into $P$.

We finish the proof using the same arguments as in the proof of Theorem \ref{theorem-totally-geodesic-sphere}.
\end{proof}

\begin{Bem}
As to be expected, similar to the case of harmonic maps, the composition $\varphi=\psi\circ\phi$ of a biharmonic map $\phi$ and a totally geodesic map $\psi$ is biharmonic.
Indeed, for arbitrary maps $\psi$ and $\phi$ and for any $\sigma\in\Gamma(\phi^\ast TN)$ we have
$$
\nabla^{\varphi}_Xd\psi(\sigma) = d\psi\left(\nabla^{\phi}_X\sigma\right) + \nabla d\psi\left(d\phi(X),\sigma\right).
$$
In particular, when $\psi$ is totally geodesic we get $\nabla^{\varphi}_Xd\psi(\sigma) = d\psi\left(\nabla^{\phi}_X\sigma\right)$, and then, by straightforward computations we find
$$
\tau_2(\psi\circ\phi)=d\psi\left(\tau_2(\phi)\right).
$$
Thus, if $\phi$ is biharmonic, $\varphi$ is biharmonic too.

In the particular case when $\psi$ defines a regular, totally geodesic submanifold, we can give an alternative proof of the above fact, similar to the proof of Theorem \ref{theorem-totally-geodesic}.
\end{Bem}

\subsection{A remark on Sampson's maximum principle}
Another geometric application of the unique continuation property of harmonic maps
was provided by Sampson \cite[Theorem 2]{MR510549}. More precisely, he gave the following result:

\begin{Theorem}[Sampson]
Assume that \(\phi\colon M\to N\) is a harmonic map, with \(q = \phi(p)\).
Let \(S\) be a \(C^2\) hypersurface in \(N\) passing through \(q\), at which point we assume that the second fundamental form is definite. If \(\phi\) is not a constant mapping,
then no neighbourhood of \(p\) is mapped entirely on the concave side of \(S\).
\end{Theorem}

Besides the unique continuation property the proof of this theorem makes use of the
maximum principle. This powerful tool only exists for second order elliptic partial differential equations, not for fourth order ones, such that we cannot expect to find a generalization of Sampson's maximum principle for proper biharmonic maps.

Indeed, in \cite{MR3045700,MR3871569}, the authors considered the map
$$
\phi_{\alpha}:\mathbb{S}^m\times \mathbb{R}\to \mathbb{R}^{m+1}, \quad \phi_{\alpha}(x,t)=\alpha(t)x
$$
and, using a reduction technique, they proved that $\phi_{\alpha}$ is biharmonic if and only if
$$
\alpha^{(4)} - 2m\alpha^{''} + m^2\alpha = 0.
$$
We note that, using the pull-back connections and the standard identifications, by direct computations we get
$$
\tau\left(\phi_{\alpha}\right) = \left( \alpha^{''} - m\alpha \right) x \quad \text{and} \quad \tau_2\left(\phi_{\alpha}\right) = \left( \alpha^{(4)} - 2m\alpha^{''} + m^2\alpha\right)x.
$$

The map $\phi_{\alpha}$ is harmonic if and only if
$$
\alpha(t) = c_1e^{\sqrt{m}t} + c_2e^{-\sqrt{m}t}, \quad t\in\mathbb{R},
$$
where $c_1$ and $c_2$ are real constants.

We note that if $\alpha(t)$ admits a local extremum point at $t_0$ such that $\alpha(t_0)>0$, then $t_0$ is a minimum point. Indeed, consider $t_0\in I$, $I$ being an open interval of $\mathbb{R}$, such that $\alpha(t)>0$ on $I$.
If $t_0$ is an extremum point, harmonicity implies $\alpha^{''}(t_0)=m\alpha(t_0)>0$ and it follows that $t_0$ is a minimum point.

The map $\phi_{\alpha}$ is biharmonic if and only if
$$
\alpha(t) = c_1e^{\sqrt{m}t} + c_2te^{\sqrt{m}t} + c_3e^{-\sqrt{m}t} + c_4te^{-\sqrt{m}t}, \quad t\in\mathbb{R},
$$
where $c_1$, $c_2$, $c_3$ and $c_4$ are real constants.

\begin{Bem}
\begin{enumerate}
\item In the biharmonic case we have solutions $\alpha(t)$ that admit a local maximum point at $t_0$ such that $\alpha(t_0)>0$. For example we may consider
$$
\alpha(t)=te^{-\sqrt{m}t} \quad \text{and} \quad t_0=\frac{1}{\sqrt{m}}.
$$
Thus, the image $\phi_{\alpha}\left(\mathbb{S}^m\times(0,\infty)\right)$ lies in the concave side of the $m$-sphere $\mathbb{S}^m\left(\frac{1}{\sqrt{m}e}\right)$ of radius $\frac{1}{\sqrt{m}e}$.
\item We also have solutions $\alpha(t)$ that admit a local minimum point at $t_0$ such that $\alpha(t_0)>0$. For example we have
$$
m=1, \quad \alpha(t)=e^t-te^{-t} \quad \text{and} \quad t_0=0.
$$
\end{enumerate}
\end{Bem}

\par\medskip
\textbf{Acknowledgements:}
The authors would like to thank Luc Lemaire and John Wood for inspiring discussions
on Sampson's maximum principle.

The first author gratefully acknowledges the support of the Austrian Science Fund (FWF)
through the project P30749-N35 ``Geometric variational problems from string theory''.

\bibliographystyle{plain}
\bibliography{mybib}
\end{document}